\documentclass[11pt]{amsart}
\usepackage{amsmath,amsthm,amsfonts,amssymb,latexsym,epsfig}

\newtheorem{theorem}{Theorem}[section]

\numberwithin{equation}{section}

\theoremstyle{definition}

\theoremstyle{remark}

\begin{document}
\title{A complement to Diananda's inequality}
\author{Peng Gao}
\address{Department of Mathematics, School of Mathematics and System Sciences, Beijing University of Aeronautics and Astronautics, P. R. China}
\email{penggao@buaa.edu.cn}
\subjclass{Primary 26D15}\keywords{Power means}
\thanks{The author is supported in part by NSFC grant 11371043.}

\begin{abstract}  Let $M_{n,r}=(\sum_{i=1}^{n}q_ix_i^r)^{\frac {1}{r}}, r \neq 0$ and $M_{n,0}=\lim_{r \rightarrow 0}M_{n,r}$ be the weighted power means of $n$ non-negative numbers $x_i$ with $q_i > 0$ satisfying $\sum^n_{i=1}q_i=1$. In particular,
   $A_n=M_{n,1}, G_n=M_{n,0}$ are the arithmetic and geometric means of these numbers, respectively.  A result of Diananda shows that
\begin{align*}
    M_{n,1/2}-qA_n-(1-q)G_n & \geq 0, \\
    M_{n,1/2}-(1-q)A_n-qG_n & \leq 0,
\end{align*}
  where $q=\min q\sb i$.
In this paper, we prove analogue inequalities in the reversed direction.
\end{abstract}

\maketitle
\section{Introduction}

   Let $M_{n,r}({\bf x}; {\bf q})$ be the weighted power means:
   $M_{n,r}({\bf x}; {\bf q})=(\sum_{i=1}^{n}q_ix_i^r)^{\frac {1}{r}}$, where
   $M_{n,0}({\bf x}; {\bf q})$ denotes the limit of $M_{n,r}({\bf x}; {\bf q})$ as
   $r\rightarrow 0$, ${\bf x}=(x_1, \ldots,
   x_n)$, ${\bf q}=(q_1, \ldots,
   q_n)$ with $x_i \geq 0, q_i>0$ for all $1 \leq i \leq n$ and $\sum_{i=1}^nq_i=1$. In this paper, unless otherwise specified, we let $q=\min q\sb i$ and
   we assume that
   $0\leq x_1 < x_2 < \cdots < x_n$.

   We define $A_n({\bf x};{\bf q})=M_{n,1}({\bf x};{\bf q}), G_n({\bf x};{\bf q})=M_{n,0}({\bf x};{\bf q}), \sigma_n=\sum_{i=1}^{n}q_i(x_i-A_n)^2$. We shall write $M_{n,r}$ for $M_{n,r}({\bf x};{\bf q})$
   and similarly for other means when there is no risk of
   confusion.

   In \cite{G5}, the following bounds of $M_{n,1/r}$ in terms of $A_n, G_n$ are given:
\begin{align}
\label{1.1}
     M_{n,\frac {1}{r}} & \leq (1-q)^{r-1} A_n +(1-(1-q)^{r-1}) G_n, \quad 1 < r \leq 2 ; \\
\label{1.2}
    M_{n,\frac {1}{r}} & \geq q^{r-1} A_n +(1-q^{r-1}) G_n, \quad r \geq 2.
\end{align}
   The reversed inequality of \eqref{1.1} is valid when $0<r<1$ and the above inequalities are generalizations of a result of Diananda (\cite{dian}, \cite{dian1}), which corresponds to case $r=2$ of the above inequalities. Note that except for the case $r=2$, the above inequalities only provide one-sided bound for any given $M_{n,1/r}$. It is therefore natural to seek for bounds that are complementary to the above ones. In this paper, we consider one way to achieve this by establishing
    the following
\begin{theorem}
\label{thm1}
   For $r \geq 2$, we have
\begin{align}
\label{1.3}
    M_{n,\frac {1}{r}}-q^{r-1} A_n -(1-q^{r-1}) G_n \leq \frac
   {1/r-q^{r-1}}{2x_1}\sigma_n,
\end{align}
   with equality holding if and only if $x_1=x_2=\cdots=x_n$ or $r=n=2,
    q=1/2$.

    For $1< r \leq 2$, we have
\begin{align}
\label{1.4}
   M_{n,\frac {1}{r}}-(1-q)^{r-1} A_n -(1-(1-q)^{r-1}) G_n \geq \frac
   {1/r-(1-q)^{r-1}}{2x_1}\sigma_n,
\end{align}
   with equality holding if and only if $x_1=x_2=\cdots=x_n$ or $r=n=2,
    q=1/2$. The reversed inequality of \eqref{1.4} holds for $1/2 \leq r<1$ with
 equality holding if and only if $x_1=x_2=\cdots=x_n$.
\end{theorem}

     Our result in fact is motivated by the the following bounds for the differences of means:
\begin{align}
\label{1.5'}
    \frac {r-s}{2x_n}\sigma_n \leq M_{n,r}-M_{n,s} \leq \frac {r-s}{2x_1} \sigma_n , \quad r>s.
\end{align}
    The above inequalities are closely related to the Ky Fan inequalities and are not valid for all $r>s$. When they are valid, then the constant $(r-s)/2$ is best possible (see \cite{G4}) and a necessary condition for inequalities \eqref{1.5'} to be valid is that $0 \leq r+s \leq 3$ (see \cite[Lemma 3.1]{G4}).
    Moreover, it is shown in \cite[Theorem 3.2]{G4} that if $r = 1$, then inequalities \eqref{1.5'} hold if and only if $-1 \leq s <1$. If $s = 1$, then inequalities \eqref{1.5'} hold if and only if $1 < r \leq 2$. In particular, the case $r=1,s=0$ of \eqref{1.5'} yields a result of Cartwright and Field \cite{C&F}:
\begin{align}
\label{1.5}
   \frac {\sigma_n}{2x_n}
    \leq   A_n-G_n  \leq  \frac {\sigma_n}{2x_1}.
\end{align}
    Using \eqref{1.5} while noting that the constant $1/2$ is best possible, one sees easily that when $r=2$, the results given in Theorem \ref{thm1} are not comparable to the bounds given by \eqref{1.1}-\eqref{1.2}.

    We can recast inequality \eqref{1.3} as
\begin{align}
\label{1.6}
    M_{n,\frac {1}{r}}-G_n-\frac
   {1/r}{2x_1}\sigma_n \leq q^{r-1} (A_n -G_n -\frac 1{2x_1}\sigma_n),
\end{align}
    from which we see that inequality \eqref{1.3} can be interpreted as a comparison between different inequalities in \eqref{1.5'}. We can deduce a similar inequality from \eqref{1.4}. This combined with our discussions above allows us to prove the right-hand side inequality of \eqref{1.5'} for $s=0, 0< r \leq 1/2$ and $1<r\leq 2$.  It is then interesting to determine all the values of $r$ such that inequalities \eqref{1.5'} hold for $r$ and $s=0$. We shall do this in Section \ref{sec 4} as we prove the following
\begin{theorem}
\label{thm2}
  Let $r \neq 0, x_1=\min \{ x_i \}, x_n=\max \{ x_i \}$, then the right-hand side inequality of \eqref{1.5'} holds with $s=0$ if and only if $0< r \leq 2$, the left-hand side inequality of \eqref{1.5'} holds with $s=0$ if and only if $1 \leq r \leq 3$. Moreover, in all these cases we have equality holding if and only if $x_1=x_2=\cdots=x_n$.
\end{theorem}

     We note that Theorem \ref{thm2} implies that the reversed inequality of \eqref{1.4} does not hold for $0<r<1/2$ in general. For otherwise we can recast it in a form similar to inequality \eqref{1.6} to deduce the validity of the right-hand side inequality of \eqref{1.5'} for $s=0, r>2$.

     We note that the following inequality
\begin{align*}
    M_{n,\frac {1}{r}}-q^{r-1} A_n -(1-q^{r-1}) G_n \geq \frac
   {1/r-q^{r-1}}{2x_n}\sigma_n,
\end{align*}
     is not valid in general as one checks easily that when $n=2, q_1=1-q, q_2=q, x_1=0, x_2=1$, the left-hand side expression above is $0$ while the right-hand side expression is not $0$ in general. Therefore, it is not possible to have a similar lower bound for the left-hand side expression in \eqref{1.3}. Similar discussions apply to \eqref{1.4} as well.

\section{Proof of Theorem \ref{thm1}}
\label{sec 3} \setcounter{equation}{0}

   Throughout this section, we assume $n \geq 2, x_1=1$ and $1<x_2< \ldots <x_n$. We will omit the discussion on the conditions for equality in each inequality as one checks easily that the desired conditions hold by going through our arguments in what follows. We first prove inequality \eqref{1.3} and we define
\begin{equation*}
  f_n({\bf x};{\bf q},q)=M_{n,\frac {1}{r}}-q^{r-1} A_n -(1-q^{r-1}) G_n - \frac
   {1/r-q^{r-1}}{2x_1}\sigma_n.
\end{equation*}
   It suffices to show $f_n({\bf x}; {\bf q}, q) \leq 0$ and we have
\begin{equation*}
  \frac {1}{q_n}\cdot\frac {\partial {f_n}}{\partial{x_n}}
= M^{1-\frac {1}{r}}_{n,\frac {1}{r}}x^{\frac {1}{r}-1}_n-q^{r-1}
-(1-q^{r-1}) G_nx^{-1}_n - (\frac 1r-q^{r-1})(x_n-A_n) :=g_n({\bf x};
{\bf q}, q).
\end{equation*}
  It suffices to show $g_n({\bf x}; {\bf q}, q) \leq 0$ as it implies $f_n({\bf x};{\bf q},q) \leq \lim_{x_n \rightarrow x_{n-1}}f_n({\bf x};{\bf q},q)$. By adjusting the value of $q$ in the expression of $\lim_{x_n \rightarrow x_{n-1}}f_n({\bf x};{\bf q},q)$ (note that it follows from \eqref{1.5} that $\frac {\partial {f_n}}{\partial{q}} \geq 0$ ) and repeating the process, it follows easily that $f_n({\bf x}; {\bf q}, q) \leq 0$.

  Now we have
\begin{align*}
\frac 1{1-q_n}\cdot\frac {\partial {g_n}}{\partial{x_n}} &=\frac {1-1/r}{1-q_n}M^{1-\frac
{2}{r}}_{n,\frac {1}{r}}x^{\frac {1}{r}-2}_n\left (q_nx^{\frac
{1}{r}}_n -M^{\frac {1}{r}}_{n,\frac {1}{r}}\right
)+(1-q^{r-1})G_nx^{-2}_n-(\frac 1r-q^{r-1}).
\end{align*}

   We make a change of variable $x_i \rightarrow y^r_i$ to recast the right-hand side expression above as
\begin{align}
\label{2.0}
 & -(1-\frac {1}{r})(q_ny_n+(1-q_n)A'_{n-1})^{r-2}A'_{n-1}y^{1-2r}_n+(1-q^{r-1}){G'}^{(1-	q_n)r}_{n-1}y^{q_nr-2r}_n-(\frac 1r-q^{r-1}) \\
\leq & -(1-\frac {1}{r})(q_ny_n+(1-q_n)A'_{n-1})^{r-2}A'_{n-1}y^{1-2r}_n+(1-q^{r-1}){A'}^{(1-	q_n)r}_{n-1}y^{q_nr-2r}_n-(\frac 1r-q^{r-1}),  \nonumber
\end{align}
   where $A'_{n-1}=A_{n-1}({\bf y'}; {\bf q'}), G'_{n-1}=G_{n-1}({\bf y'}; {\bf q'})$, and
\begin{align*}
   {\bf y'}=(y_1, \ldots, y_{n-1}), \quad {\bf q'}=(\frac {q_1}{1-q_n}, \ldots,
   \frac {q_{n-1}}{1-q_n}).
\end{align*}
   We further denote $z=y_n/A'_{n-1}$ to see that the right-hand side expression of \eqref{2.0} is
\begin{align*}
\leq & y^{-r}_n \left (  -(1-\frac {1}{r})(q_nz+1-q_n)^{r-2}z^{1-r}+(1-q^{r-1})z^{q_nr-r}-(\frac 1r-q^{r-1})z^r{A'}^{r}_{n-1} \right ) \\
\leq & y^{-r}_n \left (  -(1-\frac {1}{r})(q_nz+1-q_n)^{r-2}z^{1-r}+(1-q^{r-1})z^{q_nr-r}-(\frac 1r-q^{r-1})z^r \right ).
\end{align*}
   It suffices to show that the last expression above is non-positive for $z \geq 1$. Note first that when $q_nr \leq 1$, the last expression above equals
\begin{align*}
 &  y^{-r}_nz^{1-r} \left (  -(1-\frac {1}{r})(q_nz+1-q_n)^{r-2}+(1-q^{r-1})z^{q_nr-1}-(\frac 1r-q^{r-1})z^{2r-1} \right ) \\
\leq & y^{-r}_nz^{1-r} \left (  -(1-\frac {1}{r})(q_n+1-q_n)^{r-2}+(1-q^{r-1})-(\frac 1r-q^{r-1}) \right ) =0.
\end{align*}

  Thus we may assume $q_nr>1$ and in this case, it suffices to show that
\begin{align*}
   u(z;q_n, q)=
   -(1-\frac {1}{r})(q_nz+1-q_n)^{r-2}z^{1-q_nr}+(1-q^{r-1})-(\frac 1r-q^{r-1})z^{(2-q_n)r} \leq 0.
\end{align*}

   Now we have
\begin{align*}
   &z^{1-(2-q_n)r} \frac {\partial u}{\partial
    z} \\
   =&
   (1-\frac {1}{r})(q_nz+1-q_n)^{r-3}z^{2-2r}\left (
   \left(rq_n-1 \right )(1-q_n)z^{-1}+\left(1-(1-q_n)r \right )q_n
   \right ) \\
   &-r(2-q_n)(\frac 1r-q^{r-1}) \\
    \leq & (1-\frac {1}{r})(q_nz+1-q_n)^{r-3}z^{2-2r}\left (
   \left(rq_n-1 \right )(1-q_n)+\left(1-(1-q_n)r \right )q_n
   \right ) -r(2-q_n)(\frac 1r-q^{r-1})  \\
   =& (1-\frac {1}{r})(q_nz+1-q_n)^{r-3}z^{2-2r}\left (2q_n-1
   \right )
   -r(2-q_n)(\frac 1r-q^{r-1}).
\end{align*}
    If $2q_n-1\leq 0$, then we have $\partial u/\partial
    z  \leq 0$, as it follows from \cite[(2.1)]{G5} that $1/r-q^{r-1} \geq 0$ when $r \geq 2$. Otherwise, note
    that
\begin{align*}
    (q_nz+1-q_n)^{r-3}z^{2-2r} \leq \max \{z^{2-2r}, z^{r-3}z^{2-2r}
    \} \leq 1.
\end{align*}
   Thus we have
\begin{align*}
   z^{1-(2-q_n)r} \frac {\partial u}{\partial
    z}  &\leq (1-\frac {1}{r})\left (2q_n-1
   \right )
   -r(2-q_n)(\frac 1r-q^{r-1})  \\
   &\leq (1-\frac {1}{r})\left (2(1-q)-1
   \right )
   -r(2-(1-q))(\frac 1r-q^{r-1}) \\
   & =  (1-\frac {1}{r})\left (1-2q \right )
   -(qr+r)\left (\frac 1r-q^{r-1} \right).
\end{align*}
   It is easy to see that the last expression above is a concave up function of $q$
   for fixed $r$ and it is $ \leq 0$ when $q=0, 1/2$. Thus $\partial u/\partial
    z \leq 0$ so that $u(z; q_n, q) \leq u(1;q_n,q)=0$. This proves inequality \eqref{1.3}.


   Now, to prove inequality \eqref{1.4},  we use
   the same notations as above to see that in this case, it suffices to show $f_n({\bf x}; {\bf q}, 1-q) \geq
   0$. Again, this follows from $\frac {\partial g_n({\bf x}; {\bf q}, 1-q)}{\partial x_n} \geq
   0$. Similar to our arguments above, it is easy to see that in this case the expression \eqref{2.0} becomes
\begin{align*}
  &  -(1-\frac {1}{r})(q_ny_n+(1-q_n)A'_{n-1})^{r-2}A'_{n-1}y^{1-2r}_n+(1-(1-q)^{r-1}){G'}^{(1-	q_n)r}_{n-1}y^{q_nr-2r}_n \\
& -(\frac 1r-(1-q)^{r-1}) :=h(y_n).
\end{align*}
   It therefore remains to show that $h(y_n) \geq 0$ for $y_n \geq A'_{n-1}$.  Note first that
\begin{align*}
   h(y_n)  & \geq -(1-\frac {1}{r}){A'}^{r-1}_{n-1}y^{1-2r}_n+(1-(1-q)^{r-1}){G'}^{(1-	q_n)r}_{n-1}y^{q_nr-2r}_n  -(\frac 1r-(1-q)^{r-1}) \\
 & = y^{1-2r}_n \left ( -(1-\frac {1}{r}){A'}^{r-1}_{n-1}+(1-(1-q)^{r-1}){G'}^{(1-	q_n)r}_{n-1}y^{q_nr-1}_n  -(\frac 1r-(1-q)^{r-1})y^{2r-1}_n \right ) \\
& :=y^{1-2r} \tilde{h}(y_n).
\end{align*}
   When $q_nr-1 \geq 0$, one checks that $\tilde{h}$ is an increasing function of $y_n$ (note that in our case $(1-q)^{r-1} \geq 1/r$), hence is minimized at $y_n=A'_{n-1}$ and it is easy to see that in this case $h(y_n) \geq 0$ is equivalent to $h(A'_{n-1}) \geq 0$.

   We now consider the case $1-rq_n>0$. Note first that $\lim_{y_n \rightarrow \infty}h(y_n) \geq 0$. If $h(y_n)$ is minimized as some $y_n=y > A'_{n-1}$, then we must have $h'(y)=0$, which yields
\begin{align*}
  & \frac 1{(2-q_n)r}(1-\frac {1}{r})(q_ny+(1-q_n)A'_{n-1})^{r-3}A'_{n-1}y^{1-2r}\left ( \left (r+1 \right )q_ny+\left ( 2r-1 \right )\left ( 1-q_n\right )A'_{n-1}\right) \\
=& (1-(1-q)^{r-1}){G'}^{(1-	q_n)r}_{n-1}y^{q_nr-2r}.
\end{align*}
   This allows us to rewrite the expression for $h(y)$ as
 \begin{align}
\label{2.2}
  h(y) =&  -(1-\frac {1}{r})(q_ny+(1-q_n)A'_{n-1})^{r-2}A'_{n-1}y^{1-2r}+(1-(1-q)^{r-1}){G'}^{(1-	q_n)r}_{n-1}y^{q_nr-2r} \\
& -(\frac 1r-(1-q)^{r-1}) \nonumber\\
=& -(1-\frac {1}{r})(q_ny+(1-q_n)A'_{n-1})^{r-2}A'_{n-1}y^{1-2r} \nonumber \\
&+\frac 1{(2-q_n)r}(1-\frac {1}{r})(q_ny+(1-q_n)A'_{n-1})^{r-3}A'_{n-1}y^{1-2r} \nonumber \\
& \cdot \left ( \left (r+1 \right )q_ny+\left ( 2r-1 \right )\left ( 1-q_n\right )A'_{n-1}\right)-(\frac 1r-(1-q)^{r-1}) \nonumber \\
 =& \frac 1{(2-q_n)r}(1-\frac {1}{r})(q_ny+(1-q_n)A'_{n-1})^{r-3}A'_{n-1}y^{1-2r} \nonumber \\
 & \cdot \left ( \left (1-r+rq_n\right)q_ny- \left (1-rq_n \right) \left (1-q_n \right)A'_{n-1}\right ) -(\frac 1r-(1-q)^{r-1}). \nonumber
\end{align}
   We have
\begin{align*}
   h(y) \geq & \frac 1{(2-q_n)r}(1-\frac {1}{r})(q_ny+(1-q_n)A'_{n-1})^{r-3}A'_{n-1}y^{1-2r} \\
   & \cdot \left ( \left (1-r+rq_n\right)q_ny- \left (1-rq_n \right)\left (1-q_n \right)y \right )-(\frac 1r-(1-q)^{r-1}) \\
=& -\frac {1-2q_n}{(2-q_n)r}(1-\frac {1}{r})(q_ny+(1-q_n)A'_{n-1})^{r-3}A'_{n-1}y^{2-2r} -(\frac 1r-(1-q)^{r-1}).
\end{align*}
   If $1-2q_n \leq 0$, then $h(y) \geq 0$. When $1-2q_n>0$, we see that
\begin{align*}
   h(y) \geq & -\frac {1-2q_n}{(2-q_n)r}(1-\frac {1}{r}){A'}^{r-3}_{n-1}\cdot A'_{n-1} \cdot {A'}^{2-2r}_{n-1} -(\frac 1r-(1-q)^{r-1}) \\
=&  -\frac {1-2q_n}{(2-q_n)r}(1-\frac {1}{r}){A'}^{-r}_{n-1}-(\frac 1r-(1-q)^{r-1}) \\
\geq & -\frac {1-2q_n}{(2-q_n)r}(1-\frac {1}{r})-(\frac 1r-(1-q)^{r-1}) \\
\geq & -\frac {1-2q}{(2-q)r}(1-\frac {1}{r})-(\frac 1r-(1-q)^{r-1}).
\end{align*}
    We want to show the last expression above is non-negative. By setting $x=1-q$, we see that this is equivalent to showing that for $1/2 \leq x \leq 1$,
\begin{align}
\label{2.3}
   m(x) := x^r+x^{r-1}-\left (\frac 3r-\frac 2{r^2} \right )x-\frac 1{r^2} \geq 0.
\end{align}
    We have
\begin{align*}
   m'(x) &=rx^{r-1}+(r-1)x^{r-2}-\left (\frac 3r-\frac 2{r^2} \right ), \\
   m''(x) &=(r-1)x^{r-3}(rx+r-2).
\end{align*}
   It is easy to see that $m'(x)\geq 0$ for $1/2 \leq x \leq 1$ when $r=2$. For $1<r<2$, we see that
\begin{align*}
   m' \left(\frac {2-r}{r} \right )=\left ( \frac {2-r}{r} \right )^{r-2}-\left (\frac 3r-\frac 2{r^2} \right ).
\end{align*}
   We want to show the above expression is non-negative and we recast it as
\begin{align*}
   \left ( \frac {2-r}{r} \right )^{\frac {2-r}{3-r}}\left (\frac 3r-\frac 2{r^2} \right )^{\frac {1}{3-r}} \leq 1.
\end{align*}
   Applying the arithmetic-geometric mean inequality, we see that
\begin{align*}
   \left ( \frac {2-r}{r} \right )^{\frac {2-r}{3-r}}\left (\frac 3r-\frac 2{r^2} \right )^{\frac {1}{3-r}} \leq \left ( \frac {2-r}{r} \right )\cdot \frac {2-r}{3-r}+\left (\frac 3r-\frac 2{r^2} \right )\cdot \frac {1}{3-r} .
\end{align*}
    It therefore suffices to show the right-hand side expression above is $\leq 1$, which is equivalent to
\begin{align*}
  (r-1)(2r-1)(r-2) \leq 0.
\end{align*}
  As $1 < r < 2$, we see that the above inequality holds, hence it follows that $m'((2-r)/r) \geq 0$. As it is also easy to see that $m'(1) \geq 0$, $m'(1/2) \geq 0$, we see that $m'(x) \geq 0$ for $1/2 \leq x \leq 1$ when $1 < r < 2$. Thus, we conclude that when $1<r \leq 2$,
$m(x) \geq m(1/2) \geq 0$.

  We then conclude that in the case $1-rq_n>0$, it also suffices to show that $h(A'_{n-1}) \geq 0$, which is
\begin{align*}
  -(1-\frac {1}{r}){A'}^{-r}_{n-1}+(1-(1-q)^{r-1}){G'}^{(1-	q_n)r}_{n-1}{A'}^{q_nr-2r}_{n-1}  -(\frac 1r-(1-q)^{r-1}) \geq 0.
\end{align*}
   We recast the above inequality as
\begin{align}
\label{2.10}
 \frac {1-(1-q)^{r-1}}{1-1/r} \left (\frac {G'_{n-1}}{A'_{n-1}} \right )^{(1-	q_n)r} +\frac {(1-q)^{r-1}-\frac 1r}{1-1/r}{A'}^{r}_{n-1} \geq 1.
\end{align}
   Applying the arithmetic-geometric mean inequality, we see that
\begin{align*}
& \frac {1-(1-q)^{r-1}}{1-1/r} \left (\frac {G'_{n-1}}{A'_{n-1}} \right )^{(1-	q_n)r} +\frac {(1-q)^{r-1}-\frac 1r}{1-1/r}{A'}^{r}_{n-1} \\
\geq &  \left (\frac {G'_{n-1}}{A'_{n-1}} \right )^{(1-q_n)r^2(1-(1-q)^{r-1})/(r-1)}\cdot {A'}^{r^2((1-q)^{r-1}-1/r)/(r-1)}_{n-1}.
\end{align*}
   It follows that in order for \eqref{2.10} to be valid, it suffices to show that
\begin{align}
\label{2.12}
   G'_{n-1} \geq  {A'}^{1-((1-q)^{r-1}-1/r)/((1-q_n)(1-(1-q)^{r-1}))}_{n-1}.
\end{align}
   Note first that the above inequality holds trivially when
\begin{align*}
  1-\frac {(1-q)^{r-1}-1/r}{(1-q_n)(1-(1-q)^{r-1})} \leq 0.
\end{align*}
   Thus, we may assume the left-hand side expression above is $>0$. For any given $q_i, 1 \leq i \leq n$ (and hence $q$), we recast inequality \eqref{2.12} as
\begin{align}
\label{2.14}
   F(x_1, x_2, \ldots, x_{n-1}):={G'}^{\frac {1-q_n}{2-q_n-(1-1/r)/(1-(1-q)^{r-1})}}_{n-1} -  A'_{n-1} \geq 0.
\end{align}
   The above inequality holds trivially when $n=2$. Assuming $n \geq 3$ and let ${\bf a}=(a_1, \ldots,
    a_{n-1}) \in [1, \infty)^{n-1}$ be the point in which the absolute
    minimum of $F$ is reached. We may assume that $1=a_1 < a_2<\ldots < a_{n-1}$. As the function $x-x^r$ is decreasing for $1/2 \leq x \leq 1$, it follows easily that
\begin{align*}
   2(1-q)(1-(1-q)^{r-1}) \leq 1-2^{1-r} \leq 1-\frac 1r.
\end{align*}
  We then deduce that
\begin{align}
\label{2.16}
  \frac {q_{n-1}}{2-q_n-(1-1/r)/(1-(1-q)^{r-1})} \geq 1.
\end{align}
   It follows that
\begin{align*}
  \lim_{x_{n-1} \rightarrow \infty}F=\infty.
\end{align*}
   Thus, $a_2, \ldots, a_{n-1}$ must solve the equation
\begin{equation*}
    \nabla F=0.
\end{equation*}
   As it is easy to see that the above equation has only one root, we conclude that $n=3$ so it remains to prove \eqref{2.14} for this case. We write $x=x_2 \geq x_1=1$ to recast inequality \eqref{2.14} in this case as
\begin{align*}
  v(x) :=x^{\frac {q_2}{2-q_3-(1-1/r)/(1-(1-q)^{r-1})}} -\frac {q_2}{1-q_3}x-\frac {q_1}{1-q_3} \geq 0.
\end{align*}

   Again, by \eqref{2.16}, we see that for $x \geq 1$,
\begin{align*}
   v'(x) \geq \frac {q_2}{2-q_3-(1-1/r)/(1-(1-q)^{r-1})}-\frac {q_2}{1-q_3} \geq 0.
\end{align*}
   Thus, we have for $x \geq 1$,
\begin{align*}
  v(x) \geq v(1)=0.
\end{align*}
  This completes the proof of inequality \eqref{1.4}.

   Lastly, for the reversed inequality of \eqref{1.4} for $1/2 \leq r <1$, the proof is similar to that of \eqref{1.4}. Here it suffices to show that $h(y_n) \leq 0$ for $y_n \geq A'_{n-1}$.  Note first that in this case we always have $1-rq_n >0$. As $\lim_{y_n \rightarrow \infty}h(y_n) \leq 0$, we see that if $h(y_n)$ is minimized at some $y_n=y > A'_{n-1}$, then we must have $h'(y)=0$, which again allows us to recast $h(y)$ as the last expression in \eqref{2.2}. As $1-r+rq_n \geq 0$, we see that
\begin{align*}
    h(y) & \leq  \frac 1{(2-q_n)r}(1-\frac {1}{r})(q_ny+(1-q_n)A'_{n-1})^{r-3}A'_{n-1}y^{1-2r} \\
    & \cdot \left ( \left (1-r+rq_n\right)q_nA'_{n-1}- \left (1-rq_n \right)\left (1-q_n \right)A'_{n-1} \right ) -(\frac 1r-(1-q)^{r-1}) \\
=& -\frac {1-2q_n}{(2-q_n)r}(1-\frac {1}{r})(q_ny+(1-q_n)A'_{n-1})^{r-3}{A'}^2_{n-1}y^{1-2r} -(\frac 1r-(1-q)^{r-1}).
\end{align*}

  It follows that $h(y) \leq 0$ when $1-2q_n \leq 0$. When $1-2q_n>0$, we have (note that $r \geq 1/2$)
\begin{align*}
   h(y) \leq & -\frac {1-2q_n}{(2-q_n)r}(1-\frac {1}{r}){A'}^{r-3}_{n-1}\cdot {A'}^2_{n-1} \cdot {A'}^{1-2r}_{n-1} -(\frac 1r-(1-q)^{r-1}) \\
=&  -\frac {1-2q_n}{(2-q_n)r}(1-\frac {1}{r}){A'}^{-r}_{n-1}-(\frac 1r-(1-q)^{r-1}) \\
\leq & -\frac {1-2q_n}{(2-q_n)r}(1-\frac {1}{r})-(\frac 1r-(1-q)^{r-1}) \\
\leq & -\frac {1-2q}{(2-q)r}(1-\frac {1}{r})-(\frac 1r-(1-q)^{r-1}).
\end{align*}
   It is easy to show that the last expression above is $\leq 0$ as the function $m(x)$ defined in \eqref{2.3} is concave up for $1/2 \leq x \leq 1$ and $m(1/2) \leq 0, m(1) \leq 0$. Thus, it remains to prove $h(A'_{n-1}) \leq 0$ and we omit the argument here as it is analogue to that of the case $1 < r \leq 2$.

\section{Proof of Theorem \ref{thm2}}
\label{sec 4} \setcounter{equation}{0}

     We consider inequalities \eqref{1.5'} with $s=0$ being fixed throughout this section. Once again we omit the discussions on the conditions for equality in each inequality we shall prove.  First note that as the right-hand side inequality of \eqref{1.5'} for $0<r \leq 1/2$ and $1 \leq r \leq 2$ follows from Theorem \ref{thm1} and \eqref{1.5}, we only need to prove it for $1/2<r<1$.
     We may assume that $x_1=1<x_2<\cdots <x_n$, $q_i>0, 1 \leq i \leq n$. Let
\begin{align*}
   D(x_1, \cdots, x_n, q_1, \cdots, q_n)=M_{n,r}-G_n-\frac r2\sigma_n.
\end{align*}
    To show $D \leq 0$, it suffices to show that
\begin{align*}
   D_1(x_1, \cdots, x_n, q_1, \cdots, q_n): =\frac {\partial D}{q_n\partial x_n}=M^{1-r}_{n,r}x^{r-1}_n-\frac {G_n}{x_n}-r(x_n-A_n) \leq 0.
\end{align*}
   When $n \geq 3$, we regard $x_1=1, x_n$ as fixed and assume that $D_1$ is maximized at some point ${\bf x}'=(x'_1, \cdots, x'_n, q'_1, \cdots, q'_n)$ with $x'_1=x_1, x'_n=x_n$. Then at this point we must have
\begin{align*}
   \frac {\partial D_1}{\partial x_i}\Big |_{{\bf x}'}=0, \quad 2 \leq i \leq n-1.
\end{align*}
   Thus, the $x'_i, 2 \leq i \leq n-1$ are solutions of the equation:
\begin{align*}
   d_1(x):=(1-r)M^{1-2r}_{n,r}x^{r-1}_nx^{r-1}-\frac {G_n}{x_nx}+r=0.
\end{align*}
   It is easy to see that the above equation can have at most two different positive roots.

   On the other hand, by applying the method of Lagrange multipliers, we let
\begin{align*}
   \tilde{D}_1(x_1, \cdots, x_n, q_1, \cdots, q_n)=D_1(x_1, \cdots, x_n, q_1, \cdots, q_n)-\lambda(\sum^n_{i=1}q_i-1),
\end{align*}
   where $\lambda$ is a constant. Then at ${\bf x}'$ we must have
\begin{align*}
   \frac {\partial \tilde{D}_1}{\partial q_i}\Big |_{{\bf x}'}=0, \quad 1 \leq i \leq n.
\end{align*}
   Thus, the $x'_i, 1 \leq i \leq n$ are solutions of the equation:
\begin{align*}
   d_2(x) := \frac {1-r}{r}M^{1-2r}_{n,r}x^{r-1}_nx^{r}-\frac {G_n}{x_n}\ln x+rx-\lambda=0.
\end{align*}
   Note that $d'_2(x)=d_1(x)$ so that there is a root of $d_1(x)=0$ between any two adjacent $x_i, x_{i+1}, 1 \leq i \leq n-1$. This would imply that
   $d_1(x)=0$ has at least three different positive roots, a contradiction.

   Thus, it suffices to show $D_1 \leq 0$ for $n=2$. In this case, we let $0< q_1=q <1, q_2=1-q, x_1=x>x_2=1$ (note that we no longer assume $q=\min \{ q_1, q_2 \}$ from now on) to recast $D_1$ as
\begin{align*}
   D_1(x, q)=(qx^r+1-q)^{(1-r)/r}x^{r-1}-x^{q-1}-r(1-q)(x-1).
\end{align*}
    Note that
\begin{align*}
   D_2(x,q) :=(1-q)^{-1}D'_1(x, q)=(r-1)(q+(1-q)x^{-r})^{(1-2r)/r}x^{-r-1}+x^{q-2}-r.
\end{align*}
   As $r-1<0$, it follows from the arithmetic-geometric mean inequality with non-positive weights that
\begin{align*}
   \frac {r-1}{r}(q+(1-q)x^{-r})^{(1-2r)/r}x^{-r-1}+\frac 1rx^{q-2}\leq (q+(1-q)x^{-r})^{(1-2r)(r-1)/r^2}x^{(q-1-r^2)/r} \leq 1.
\end{align*}
   This implies that $D'_1(x, q) \leq 0$ and hence $D_1(x,q) \leq 0$ for $x \geq 1$ and this completes the proof for the right-hand side inequality of \eqref{1.5'} for $1/2<r<1$.

   Next, note that
\begin{align*}
  \lim_{q \rightarrow 0^+}\frac {D(x,1, q, 1-q)}{q}=\frac {x^r-1}{r}-\ln x-\frac {r}{2}(x-1)^2 .
\end{align*}
    As the right-hand side expression above is positive when $r>2$ and $x \rightarrow +\infty$ , we then conclude that in order for the right-hand side inequality of \eqref{1.5'} for $s=0$ to hold, it is necessary to have $r \leq 2$ and this completes the proof for the assertion on the right-hand side inequality of \eqref{1.5'} for $s=0$.

    Note we also have
\begin{align*}
  \lim_{q \rightarrow 1^-}\frac {D(x,1, q, 1-q)}{1-q}=x\ln x-\frac {x-x^{1-r}}{r}-\frac {r}{2}(x-1)^2.
\end{align*}
    As the right-hand side expression above is negative when $0<r<1$ and $x \rightarrow 0^+$, we conclude that in order for the left-hand side inequality of \eqref{1.5'} for $s=0$ to hold, we must have $r \geq 1$. As \cite[Lemma 3.1]{G4} implies that we also need to have $r \leq 3$ in this case, this shows that it is necessary to have $1 \leq r \leq 3$ in order for the left-hand side inequality of \eqref{1.5'} for $s=0$ to hold.

  It remains to prove the left-hand side inequality of \eqref{1.5'} for $s=0, 1 \leq r \leq 3$. Note that the case $1 \leq r \leq 2$ is a consequence of the left-hand side inequality of \eqref{1.5} and the left-hand side inequality of \eqref{1.5'} for $s=1, 1 < r \leq 2$, valid according to \cite[Theorem 3.2]{G4}. Thus, we may assume that $2 < r \leq 3$. In this case, it suffices to show $D \geq 0$ provided that we assume $0<x_1<x_2<\cdots <x_n=1$. Similar to our discussions above, one shows easily that this follows from $\partial D/\partial x_1 \leq 0$ for $n=2$, which is equivalent to $D_1(x,q) \leq 0$ for $0<x \leq 1$. As $D_1(1,q)=0$, it suffices to show that $D_2(x,q) \geq 0$ for $0<x \leq 1$. As $\lim_{x \rightarrow 0^+}D_2(x,q)>0, D_2(1,q)=0$, we only need to show the values of $D_2$ at points satisfying:
\begin{align*}
   \frac {\partial D_2}{\partial x}=0,
\end{align*}
  are non-negative.

   Calculation shows that at these points, we have
\begin{align*}
   (r-1)(qx^r+1-q)^{(1-2r)/r}x^{r-2}=\frac {(qx^r+1-q)(2-q)x^{q-2}}{-q(r+1)x^r+(r-2)(1-q)}.
\end{align*}
   We may assume the denominator of the right-hand side expression above is positive. Substituting this back to the expression for $D_2(x,q)$, we see that it remains to show that for $0 < x \leq 1$,
\begin{align*}
   \frac {(qx^r+1-q)(2-q)x^{q-2}}{-q(r+1)x^r+(r-2)(1-q)}+x^{q-2}-r \geq 0.
\end{align*}
   As $x^{q-2} \geq 1$ for $0 < x \leq 1$, it suffices to show for $0 < x \leq 1$,
\begin{align*}
   \frac {(qx^r+1-q)(2-q)x^{q-2}}{-q(r+1)x^r+(r-2)(1-q)}\geq r-1.
\end{align*}
   We recast the above inequality as
\begin{align*}
   j(x,q):=q(2-q)x^{r+q-2}+(1-q)(2-q)x^{q-2}+q(r^2-1)x^r -(1-q)(r-1)(r-2)\geq 0.
\end{align*}
   Again as $\lim_{x \rightarrow 0^+}j(x,q)>0, j(1,q) \geq 0$ when $2 < r \leq 3$, we only need to show the values of $j(x,q)$ at points satisfying:
\begin{align*}
   \frac {\partial j}{\partial x}=0,
\end{align*}
  are non-negative.

  Calculation shows that at these points, we have
\begin{align*}
   q(r^2-1)x^r=\frac {(1-q)(2-q)^2x^{q-2}}{r}-\frac {q(2-q)(r+q-2)x^{r+q-2}}{r} .
\end{align*}
   Substituting this back to the expression for $j(x,q)$, we see that it suffices to show that
\begin{align*}
  j_1(x,q):=(1-q)(2-q)(r+2-q)x^{q-2}+q(2-q)^2x^{r+q-2} -(1-q)r(r-1)(r-2) \geq 0.
\end{align*}
   One checks easily that $\lim_{x \rightarrow 0^+}j_1(x,q) \geq 0$ and that on setting $y=1-q$, we have
\begin{align*}
  j_1(1,q)&=(1-q)(2-q)(r+2-q)+q(2-q)^2-(1-q)r(r-1)(r-2) \\
          &=1+(r+2-r(r-1)(r-2))y+(r+1)y^2 \geq 1-y+(r+1)y^2 \geq 0,
\end{align*}
   as one checks easily that $r+2-r(r-1)(r-2)$ is a decreasing function of $2<r \leq 3$, hence is minimized at $r=3$.

    Thus, we only need to show the values of $j_1(x,q)$ at points satisfying:
\begin{align*}
   \frac {\partial j_1}{\partial x}=0,
\end{align*}
   are non-negative.

   Calculation shows that at these points, we have
\begin{align*}
   q(2-q)^2x^{r+q-2}=\frac {(1-q)(2-q)^2(r+2-q)x^{q-2}}{r+q-2}.
\end{align*}
   Substituting this back to the expression for $j_1(x,q)$, we see that it suffices to show that for $0<x \leq 1$,
\begin{align*}
  \frac {(2-q)(r+2-q)x^{q-2}}{r+q-2} \geq (r-1)(r-2).
\end{align*}
   The above inequality is valid as one checks easily that when $2<r \leq 3$,
\begin{align*}
  \frac {(2-q)(r+2-q)x^{q-2}}{r+q-2} \geq \frac {(2-q)(r+2-q)}{r+q-2} \geq  \frac {r+1}{r-1}\geq (r-1)(r-2).
\end{align*}
  We now conclude that the left-hand side inequality of \eqref{1.5'} is valid for $s=0, 2<r \leq 3$ and this completes the proof of the theorem.



\end{document}